\documentclass[a4paper, 11pt]{article}
\usepackage{authblk}
\usepackage{amsmath,amssymb, amsfonts, textcomp, eufrak}
\usepackage[titletoc,toc,title]{appendix}
\newtheorem{theorem}{Theorem}[section]
\newtheorem{definition}[theorem]{Definition}
\newtheorem{lem}[theorem]{Lemma}
\newtheorem{prop}[theorem]{Proposition}
\newtheorem{rem}[theorem]{Remark}
\newtheorem{cor}[theorem]{Corollary}
\newtheorem{exa}[theorem]{Example}

\numberwithin{equation}{section}
\usepackage{graphicx}
\numberwithin{equation}{section}
\usepackage{float} 
\usepackage{hyperref}
\usepackage{fancyhdr}
\usepackage{caption}
\usepackage{subcaption}
\usepackage{makeidx}
\makeindex
\newcommand{\R}{\mathbb{R}}

\newcommand{\N}{\mathbb{N}}

\newcommand{\Rp}{\mathbb{R}^+}
\newcommand{\F}{\mathcal{F}}

\newcommand{\PS}{(\Omega, \mathcal{F}, \mathbb{P})}

\newcommand{\I}{\mathbb{I}}

\newcommand{\Bb}{\mathcal{B}}

\newcommand{\Z}{\mathbb{Z}}

\newcommand{\T}{\mathbb{T}}

\newcommand{\p}{\mathbb{P}}

\newcommand{\bn}{\begin{definition}}
\newcommand{\en}{\end{definition}} 
\newcommand{\bt}{\begin{theorem}}                
\newcommand{\et}{\end{theorem}}
 \newcommand{\bnm}{\begin{enumerate}}              
\newcommand{\enm}{\end{enumerate}}
\newcommand{\br}{\begin{rem}} 
\newcommand{\er}{\end{rem}}

\newcommand{\om}{\omega}

\newcommand{\Om}{\Omega}

\newcommand{\btm}{\begin{itemize}}
 \newcommand{\etm }{\end{itemize}}

\newcommand{\E}{\mathbb{E}}

\newcommand{\Rd}{\R^d}

\newcommand{\qed}{\hfill $\square$}
\begin{document}
\title{Random periodic solutions and ergodicity for stochastic differential equations}
\author[1]{Kenneth Uda}
\author[2]{Huaizhong Zhao}
 \affil[1]{School of Mathematics, University of Edinburgh, UK}
 \affil[2]{Department of Mathematical Sciences, Loughborough University, UK}
 
\providecommand{\keywords}[1]{\textbf{\textit{Keywords:}} #1}
\date{}
\maketitle
\begin{abstract}\rm In this paper, we establish some sufficient conditions for the existence of stable random periodic solutions of stochastic differential equations on $\Rd$ and ergodicity in the random periodic regime. The techniques involve the existence of Lyapunov type function, using two-point generator of the stochastic flow map, strong Feller argument and weak convergence. 
\end{abstract}
\newcommand{\Le}{\rho_0^{\varepsilon}}
%\tableofcontents
\section{Introduction}
Random dynamical system (RDS) is a dynamical system with some randomness, its idea was discussed in 1945 by Ulam and von Nuemann \cite{Ulam}  and few years later by Kakutani \cite{Kak} and continued in the 1970s  in the framework of ergodic theory. Series of works by Elworthy, Meyer, Baxendale, Bismut, Ikeda, Kunita, Watanabe and others \cite{Bax, Bism, Elw2, Elw1, Ikeda, Hkunita, Kunita} showing that stochastic differential equations (SDEs) induce stochastic flows, gave a substantial push to the subject. Towards late 1980s, it became clear that the techniques from dynamical systems and probability theory could produce a theory of RDS. It was extensively developed by Arnold \cite{Arnold} and his "Bremen group" based on two-parameter stochastic flows generated by stochastic differential equations due to Kunita \cite{Hkunita, Kunita} and others. 

To investigate the long time behaviour of RDS is of great interests both in applications and theory. There are two main issues that motivate the behaviour of a mathematical model in the long run with theoretical and practical consequences. One is to understand the random equilibrium and their distributions, to describe invariant property under the transformation of RDS and where the orbits (ensemble of trajectories) converge to, in the long run. Another one is to ascertain whether the limiting behaviour is still essentially the same after small changes to the evolution rule. 

 Intuitively, the limiting behaviour of a dynamical system is captured by the concept of stationary, periodic solutions or more general, invariant or quasi-invariant manifolds of  stationary or periodic solutions. %For a dynamical system $\Phi:\T\times \mathbb{M}\rightarrow \mathbb{M},$ over $t\in\T,$ a stationary solution is a point $x\in X,$ such that \begin{equation}
  % \Phi(t,x) = x, \quad \text{for all $t$}\in\T.
  % \end{equation}
%And a periodic solution is a periodic function $u:\T\rightarrow X$ with period $\tau\neq 0,$ such that \begin{equation}
%u(t+\tau) = u(t), \quad \text{and}\quad \Phi(t,u(s))= u(t+s), \quad \text{for all $t,s$}\in\T.
%\end{equation} 
To understand and give the existence of such solutions attracted vast interest in theory and applications. Periodic solutions have been crucial in the qualitative theory of dynamical systems and its systematic consideration was initiated by Poincar\'e in his work \cite{Ponc}. Periodic solutions have been studied for many fascinating physical problems, for example, van der Pol equations \cite{Van}, Li\'enard equations \cite{Lien}, etc. However, once the influence of noise on the system is considered, which is evidently inevitable in many situations; the dynamics start to depend on both time and the noise path, so stationary and periodic solutions in the usual sense may not exist for randomly perturbed systems. 

As in the deterministic systems setting, random stationary solutions are central in the long time behaviour of RDS. For an RDS $\Phi: \T^+\times\Om\times \mathbb{M}\rightarrow \mathbb{M},$ over a metric dynamical system $(\Om, \F,\p,(\theta_t)_{t\in\T}),$ a random stationary solution is an $\F$-measurable random variable $Y:\Om\rightarrow \mathbb{M}$ such that \begin{equation}\label{Stat}
\Phi(t,\om,Y(\om)) = Y(\theta_t\om), \quad \text{for all $t$}\in\T^+, \quad \p-\text{a.s.}
\end{equation}
Here $\T$ is a two-sided time domain (discrete or continuous) and $\T^+= \T\cap[0, +\infty).$
The notion of random stationary solutions of RDS is a natural extension of fixed point solution of deterministic systems. It is a one-force, one-solution setting that describes the pathwise invariance of the system over time along the base dynamical system $\theta$ on noise space and the pathwise limts of  random dynamical systems (e.g., \cite{Caraballo, Schmalfuss, Huaiz}).

Analogous to dynamical systems, the notion of random periodic solutions plays a similar role to RDS. In the physical world around us (e.g.,~biology, chemical reactions, climatic dynamics, finance, etc.), we encounter many phenomena which repeat after certain interval of time.  Due to the unavoidable random influences, many of these phenomena may be best described by random periodic paths rather than periodic solutions. For example, the maximum daily temperature in any particular region is a random process, however, it certainly has periodic nature driven by divine clock due to the rotation of the earth around the sun. There had been few attempts in physics to study random perturbation of limit cycle for some time (e.g.,~\cite{Kurr, Jeffrey}). One of the challanges that hindered real progress was the lack of a rigorous mathematical definition of random periodic solution and appropriate mathematical tools. For a random path with some periodic property, it is not obvious what a reasonable mathematical relation between the random positions $S(s,\om)$ at time $s$ and $S(s+\tau,\om)$ at time $s+\tau$ after a period $\tau$ should be. However, as $S(t,\om)$ is a true path, so it is not necessarily true that $S(s,\om) = S(s+\tau,\om).$ To require that $S(s+\tau,\om)$ is in the neighbourhood of $S(s,\om)$ by considering a small noise perturbation was worthwhile attempt. However, this approach does not apply to many stochastic differential equations and lack rigour. It is not even true for small noise and the scope of applications is limited.  Recently, in Zhao and Zheng \cite{Zhao}, Feng, Zhao and Zhou \cite{Chung}, Feng and Zhao \cite{Feng}, it has been observed that for fixed $s,$ $(S(s+k\tau,\om))_{k\in\Z}$ should be a random stationary solution of the discrete RDS $\Phi(k\tau, \om)$. This then led to the rigorous definition of random periodicity $S(s+\tau,\om)= S(s,\theta_{\tau}\om).$  For an RDS $\Phi$ over a metric dynamical systems $(\Om,\F,\p,(\theta_s)_{s\in\T}),$ a \textbf{random periodic solution}  is an $\F$-measurable function $S: \T\times\Om\rightarrow \mathbb{M},$ of period $\tau$ such that \begin{equation}\label{intreq}
S(s+\tau,\om)= S(s,\theta_{\tau}\om) \quad \text{and}\quad \Phi(t,\theta_s\om, S(s,\om)) = S(t+s, \om), \quad \text{for all $s$}\in\T, \; t\in \T^+.
\end{equation}
The study of random periodic solutions is more fascinating and difficult than deterministic periodic solutions. An extra essential difficulty is from the fact that trajectory (solution path) of the random dynamical systems does not follow a periodic path, but the pullback path $\{\eta(s,\om):=S(s,\theta_{-s}\om): 0\leqslant s\leqslant\tau\}$ is a periodic curve and random periodic path moves from one periodic curve to another one corresponding to different $\om.$ If one considers a family of tajectories starting from different points on the closed curve $\eta(.,\om),$ then the whole family of trajectories at time $t\in\T$ will lie on a closed curve $\eta(.,\theta_{\tau}\om).$ 

Existence of random periodic solutions of stochastic (partial) differential equations were investigated by Feng, Zhou and Zhao \cite{Chung} and Feng and Zhao \cite{Feng}. Their results are based on infinite horizon stochastic  integral equations and an Wiener-Sobolev compact embedding argument. In fact, one of the technical assumptions in their works was some boundedness conditions on the vector fields associated with the stochastic (partial) differential equations, though these conditions can be removed given more conditions such as weak dissipativity \cite{Zhao2}.  We employ the Lyapunov function technique to characterize the boundedness conditions (dissipativity of the stochastic flow) to prove the existence of unique stable random periodic solution. The conditions of our results are quite natural to some applicable SDEs and verifiable in terms of their coefficients. It is proved in a number of works (e.g., \cite{Mao} and references therein) that the technique of Lyapunov function could be used to provide a bound to the top Lyapunov exponent of stochastic flows. We previously employed semiuniform ergodic theory approach to prove the existence of random periodic solutions in \cite{Uda} where the conditions on the bound of top Lyapunov exponent was central.~The Lyapunov function technique discussed here makes it easier to establish ergodicity of periodic measures induced by the random periodic solutions.

The purpose of ergodicity is to study invariant measures and related problems. It is one of the well studied problems in dynamical systems, stochastic analysis, statistical physics and related areas. Roughly speaking, an ergodic dynamical system is one that its behaviour averaged over long time is the same as its behaviour over phase space. In the framework of RDS, ergodicity is the cornerstone in the investigation of long time behaviour, various techniques and their variants have been developed by many researchers (e.g.,\cite{Uarnold, Arnold, Bax2, Bax3, Hans3, Hasm, Mattingly, Meyn93, Meyn931} and references therein). Many important results were established in the regime of (random) stationary measures and (random) stationary process. However, various assumptions involved, automatically exclude random periodic regime. Feng and Zhao \cite{Zhao2} defined periodic measures and proved the Krylov-Bogolyubov procedure as a variant of Poincar\'e Bendixson theorem in the RDS framework. The Krylov-Bogolyubov procedure for periodic measures was also investigated by Hasminksii \cite{Hasm}, in terms of transition probability function. We shall recover the ergodicity of periodic measures in \cite{Zhao2} for some SDEs using the Lyapunov function technique.
 
The outline of the paper is as follows. In Section \ref{Sec2}, we present some standard notation and definitions that will be employed in our proofs. In Seection \ref{Sec3}, we prove the existence of random periodic solutions by employing the two-point generator and Lyapunov function techniques. Section \ref{Sec4} is about periodic measures induced by the random periodic solutions and their ergodicity in a certain Poinc\'are section.

\section{ Preliminaries }\label{Sec2}
In this section, we fix notation that will be frequently used throughout this paper. We also, introduce random periodic solutions for stochastic flows generated by time dependent SDEs with a simple example to fix the idea (for more examples, see \cite{Chung, Feng, Zhao2, Uda}). 
%\subsection{Notation}

On a complete separable metric space $(\mathbb{M}, \text{d})$, we denote
 the set of bounded measurable real-valued functions by
%\begin{align*} 
 $\Bb_b(\mathbb{M})$
 % =\left\{f: \mathbb{M}\rightarrow\R, \; \text{measurable}: \Vert f\Vert_{\infty}<\infty\right\},\qquad
 with the norm 
$\Vert f\Vert_{\infty} := \sup_{x\in \mathbb{M}}\vert f(x)\vert$;
%\end{align*}
and the set of bounded continuous real-valued functions by 
%\begin{align*} 
 $\mathcal{C}_b(\mathbb{M})$;
 Let $\mathcal{C}_b^{l}(\Rd;\Rd)$ be the Banach space of the functions $f: \Rd\rightarrow \Rd$ which has l-th derivative being continuous
 with the norm 
\begin{align*}
&\Vert f\Vert_{l}:= \sup_{x\in\Rd}\frac{\vert f(x)\vert}{1+\vert x\vert}+\sum_{1\leq \vert \alpha\vert \leq l}\sup_{x\in \Rd}\vert D^{\alpha}f(x)\vert.
\end{align*}
Let $(\mathbb{M}, \Bb(\mathbb{M}), \mu)$ be a Borel measure space, we denote $L^p(\mathbb{M}), \; 1\leqslant p<\infty$ as the set of real-valued Lebesgue integrable functions 
%$f: \mathbb{M}\rightarrow \R$ defined by 
with the norm
%\begin{align*}
%L^p(\mathbb{M}):=\{f\in \Bb_b(\mathbb{M}): \Vert f\Vert_{p}<\infty\}, \qquad 
$\Vert f\Vert_{p}=\left(\int_{\mathbb{M}}\vert f\vert^pd\mu\right)^{1/p}$.
%\end{align*}

In what follows, we will consider the case $\mathbb{M}= \Rd$ from time to time without causing confusions.
Let $\PS$ be a complete probability space and $\mathcal{G}\subseteq\F,$ we denote $L^p(\Om, \mathcal{G}, \p),\;p\geqslant 1$ as the space of $\mathcal{G}$-measurable random variables $X:\Om\rightarrow\R^d, \; d\in\N,$ such that $\E\vert X\vert^p<\infty,$ equiped with the $L^p$ norm 
$\Vert X\Vert_p = \left(\E\vert X\vert^p\right)^{1/p}$. 

We shall fix the probability space $\PS$ as the classical Wiener space, i.e., $\Om=\mathcal{C}_0(\R;\R^m), \; \break m\in \N$, is a linear subspace of continuous functions that take zero at $t=0$, endowed with compact open topology defined via
\begin{align*}
\text{d}(\om,\hat{\om})= \sum_{n=0}^\infty\frac{1}{2^n}\frac{\Vert \om -\hat{\om}\Vert_n}{1+\Vert \om -\hat{\om}\Vert_n}, \qquad \Vert \om -\hat{\om}\Vert_{n} = \sup_{t\in [-n, n]}\vert \om(t) -\hat{\om}(t)\vert.
\end{align*}
The sigma algebra $\F$ is the Borel sigma algebra generated by open subsets of $\Om$ and $\p$ is the Wiener measure, i.e., that the law of the process $\om\in\Om$ with $\om(0) =0.$
Let the $(\Bb(\R)\otimes\F, \F)$-measurable flow $\theta: \R\times\Om\rightarrow\Om,$  be defined by
\begin{align}\label{Mes_fl}
\theta_t\om(\cdot) = \om(t+\cdot)- \om(t).
\end{align}
It is well known that the measurable flow $(\theta_t)_{t\in \R}$ is ergodic (e.g., \cite{Arnold, Arnoldp}). Let $\F_s^t := \sigma(\om(u)-\om(v): s\leq u, v\leq t)\cup \mathcal{N},$ where $\mathcal{N}$ is a collection of $\p$-null sets of $\F,$ then
 $\F^t_s$ is a two parameter filtration on $(\Om, \F,\p)$ and  the metric dynamical system $\theta = (\Om, \F, \p, (\F^t_s)_{t\geqslant s},(\theta_t)_{t\in\R})$  is a filtered dynamical system. Note  $\theta^{-1}_{r}(\F^t_s) = \F^{t+r}_{s+r}, \; r\in \R, \quad -\infty<s\leq t<\infty.$

We consider time dependent SDE on $\Rd$ of the form
\begin{equation}\label{phnl}
dX= f_0(t,X)dt+\sum_{k=1}^mf_k(t,X)dW_t^k,\quad X(t_0)=x\in\Rd, \quad t\geqslant t_0.
\end{equation}
\begin{prop}[Stochastic flows \cite{Hkunita}]\;\label{flow property}\rm
Suppose that the coefficients $f_k(t,x),\; 0\leqslant k\leqslant m,$  of SDE (\ref{phnl}) are continuous in $t$ and uniformly Lipschitz continuous with respect to $x\in\Rd.$ Then there exists a modification of the solution of SDE (\ref{phnl}), denoted  by $X(t,t_0,\om,x)$ which satisfies the following properties:
\btm
\item[(1)] For each $t, t_0\in \T\subseteq\R, \; t\geqslant t_0$ and $x,$ $X(t,t_0,.,x)$ is $\F^{t}_{t_0}$-measurable,
\item[(2)] for almost all $\om,$ $X(t,t_0,\om,x)$ is continuous in $(t,t_0,x)$ and satisfies $$X(t_0,t_0,\om,x) = x,$$
\item[(3)] for almost all $\om,$ \begin{equation}
X(t+u,t_0,\om,x)= X(t+u,t,\om,X(t,t_0,\om,x))
\end{equation} 
is satisfied for all $t,t_0\in\T\subseteq\R, \; t\geqslant t_0$ and $u>0,$
\item[(4)] for almost all $\om,$ $X(t,t_0,\om,.):\R^d\rightarrow\R^d$ is a homeomorphism for all $t\geqslant t_0$,
\item[(5)] if in addition, the coefficients $f_0, f_1, \cdots, f_m$ are differentiable in $x$ and their first derivatives are continuous and bounded with respect to $(t,x)$, then, for almost all $\om,$ $X(t,t_0,\om,.): \Rd\rightarrow\Rd$ is a diffeomorphism for all $t\geqslant t_0.$
\etm
\end{prop}

\begin{definition}[Random periodic solution for stochastic flows \cite{Chung, Feng, Zhao}]\;\label{Chu}\rm
 A random periodic solution of period $\tau$ of a stochastic flow $X:\Delta\times\Om\times \Rd\rightarrow \Rd$ is an $\F$-measurable function $S:\T\times\Om\rightarrow\Rd$ such that $$ S(t+\tau,\om) = S(t,\theta_{\tau}\om)\quad \text{and} \quad X(t+s,s,\om,S(t,\om)) = S(t+s,\om),\quad a.s.,$$ for any $t,s \in\T; t\geqslant s$ and $\om\in\Om,$ where $\Delta:=\{(t,s)\in\T^2; t\geqslant s\}.$
\end{definition}

\begin{rem}\rm\;
 Suppose that $\{X(t,s,\om,.); t\geqslant s\}$ is a {\it time homogeneous flow}, for example, a solution of an autonomous SDEs driven by indpepndent Brownian motions. In this case, we can write \begin{align*}X(t,s,\om,x)= X(t-s,0,\theta_{s}\om,x) \quad\text{and}\quad
  \Phi(t,\om,x):= X(t,0,\om,x).
 \end{align*}
 \btm
 \item[(i)] Let $S(s,\om)$ be a random periodic solution with period $\tau$ in the sense of definition \ref{Chu}, we have \begin{align}
\nonumber \Phi(t,\theta_s\om,S(s,\om))= X(t,0,\theta_s\om,S(s,\om))&=X(t+s-s,0,\theta_s\om,S(s,\om))\\&=X(t+s,s,\om,S(s,\om))= S(t+s,\om),
 \end{align}
 corresponding to the definition of random periodic solution we have in the introduction (equation (\ref{intreq})) for a cocycle. 
 \item[(ii)] On the other hand, if we set $\eta(s,\om):=S(s,\theta_{-s}\om),$ so that, \begin{align*}\eta(s+\tau,\om)&=S(s+\tau,\theta_{-s-\tau}\om)=S(t,\theta_{\tau}\circ\theta_{-t-\tau}\om)= \eta(t,\om),
 \end{align*}and using the fact that $S(s,\om)$ is a random periodic solution of a homogeneous flow $X(t,s,\om,.)$, we have that 
 \begin{align*} 
 \eta(t+s,\theta_t\om)=S(t+s,\theta_{-s}\om)=X(t+s,s,\theta_{-s}\om,S(s,\theta_{-s}\om))\\=X(t,0,\om,\eta(s,\om))
 =\Phi(t,\om,\eta(s,\om)).
 \end{align*} 
 Thus,
 \begin{equation}
 \eta(s+\tau,\om)= \eta(s,\om) \quad \text{and}\quad \Phi(t,\om,\eta(s,\om))= \eta(t+s,\theta_s\om),
 \end{equation}
 corresponding to random periodicity in the pullback sense considered in \cite{Zhao} and in the more recent works \cite{Uda, Wan2, Wan3}.
 \etm
\end{rem}

\begin{exa}\label{ex32}\;\rm
 Consider the following SDE 
\begin{equation}\label{ehnl}
 dX = -\alpha(t)Xdt + dW_t, \quad X(t_0)= x_0\in \Rd, \quad t\geqslant t_0,
\end{equation}
 where $\alpha: \R\rightarrow \R$ is a continuous function and there exists $\tau> 0 $ such that $\alpha(t+\tau) = \alpha(t)$ with $$\int_{-\infty}^{t}e^{-2\int_{s}^{t}\alpha(u)du}ds<\infty, \quad \text{for}\quad 0\leqslant t\leqslant \tau.$$ The random variable $S(s,\om)$ defined by $$ S(s,\om) = \int_{-\infty}^{s}e^{-\int_{r}^{s}\alpha(u)du}dW_r(\om)$$ is a random periodic solution of the stochastic flow $X(t,t_0,\om,x_0)$ generated by SDE (\ref{ehnl}) defined by $$X(t,t_0,\om,x_0) = x_0e^{-\int_{t_0}^{t}\alpha(u)du} + \int_{t_0}^{t}e^{-\int_{s}^t\alpha(u)du}dW_s(\om).$$   
Indeed, by suitable change of variable and the periodicity of $\alpha,$ we have that  \begin{align*} S(s,\theta_{\tau}\om) = \int_{-\infty}^{s}e^{-\int_{r}^{s}\alpha(u)du}dW_{r+\tau}= \int_{-\infty}^{s+\tau}e^{-\int_{r}^{s+\tau}\alpha(u)du}dW_r= S(s+\tau,\om),
 \end{align*}
 and
\begin{align*} 
 X(t+s,s,\om,S(s,\om))&= e^{-\int_{s}^{t+s}\alpha(u)du}\int_{-\infty}^{s}e^{-\int_{r}^{s}\alpha(u)du}dW_r(\om)+ \int_{s}^{t+s}e^{-\int_{r}^{t+s}\alpha(u)du}dW_r(\om)\\ &=\int_{-\infty}^{t+s}e^{-\int_{r}^{t+s}\alpha(u)du}dW_r(\om) = S(t+s,\om). \quad \quad \quad \square
\end{align*}

\end{exa}

\section{Existence of random periodic solutions}\label{Sec3}
We adopt the approach of studying the infinitesimal separation of trajectories of SDEs via their Markov evolution to prove the existence of stable random periodic solutions. For this, we recall a standard notion of the transition probability function $P(t_0, x; t, A)$ induced by solutions of SDE (\ref{phnl}),
 \begin{equation}
 P(t_0,x; t, A) = \p\big(\{\om\in\Om: X(t, t_0,\om,x )\in A\}\big), \quad t\geqslant t_0, \quad A\in\Bb(\Rd).
 \end{equation}
 %The transition probability function  $P(t_0, x; t, A)$ is the probability that the process $X(t,t_0,\om,x)$ takes value in $A$ at time $t$ under the condition that it is at the point $x$ at time $t_0\leq t.$
 The Markov evolution $T_{t,t_0}: \mathcal{C}_b(\Rd)\rightarrow\mathcal{C}_b(\Rd)$ is given by 
 \begin{align*}
 T_{t, t_0}h(x) = \int_{\Rd}h(y)P(t_0,x; t, dy) = \E [h(X(t,t_0,\om,x))]
 \end{align*}
 and for any probability measure $\mu$ on $(\Rd,\Bb(\Rd)),
 $ \begin{equation}
 (T_{t,,t_0}^*\mu)(A) = \int_{\R^2}P(t_0,x; t, A)\mu(dx), \quad \text{for any $t\geqslant t_0$ and  $B$}\in\Bb(\Rd).
 \end{equation}
 
 \bigskip
 
The poineering work of Has'minskii \cite{Hasm} championed the stability theory of SDEs, by systematically adopting the concept {\it Lyapunov function} $V$ for the SDEs. The flavour in this concept is the fact that the average growth of a function $V$ along the trajectory $X(t, t_0, \om,x)$ is expressed by \begin{equation}\label{one point}
 \mathcal{L}V(t_0,x) = \lim_{t\downarrow t_0}\frac{\E[V(t, X(t,t_0, \om, x))]-V(t_0, x)}{\vert t-t_0\vert}.
 \end{equation}
For $V\in \mathcal{C}_b^{1,2}(\R\times\Rd),$ we can use It\^o's formula to write $\mathcal{L}$ as 
$$\mathcal{L}V(t,x) =\frac{\partial V(t,x)}{\partial t}+ \sum_{i=1}^{d}f_0^i(t, x)\frac{\partial V(t,x)}{\partial x_i}+ \frac{1}{2}\sum_{i,j=1}^{d}\sum_{k=1}^mf_k^i(t,x)f_k^j(t,x)\frac{\partial^2V(t,x)}{\partial x_i\partial x_j}.$$  
The generator $\mathcal{L}$ determines the law of the one-point motions $\{X(t,t_0,\om,x): x\in \Rd,\;t\geqslant t_0\}.$ 
Has\'minskii \cite{Hasm} established the stability of solution of SDEs using the differential operator $\mathcal{L}$, more work in this direction can be found in (e.g., \cite{Mao}). The idea in Has\'minskii's technique is to investigate the difference between two trajectories of a stochastic differential equation, where one of the trajectories is a deterministic (zero) solution. However, due to random fluctuations, zero solution of an SDE may not exist in general, for instance, the simple SDE in example \ref{ex32}. In order to investigate the infinitesimal separation of two nontrivial trajectories, one requires a differential operator that  gives information about the joint distribution between these trajectories. From the definition of $\mathcal{L},$ the term $f_k^i(t,x)f_k^j(t,x)$ in the sum only contains the diagonal entries of the inifinitesimal covariance of the stochastic flow $\{X(t,t_0,\om,.): t\geqslant t_0\},$ so $\mathcal{L}$ does not determine the law of stochastic flows (cf.~\cite{Arnold, Bax2, Kunita}). 

It is known (e.g., \cite{Arnold, Bax2, Kunita}) that the law of stochastic flows driven by Brownian motion is determined by the generator $\mathcal{L}^{(2)}$ of the two-point motions $\{(X(t,t_0,\om, x), X(t,t_0,\om,y)):  (x, y) \in \R^{2d}, t\geqslant t_0\}$.~We denote the transition probability function and the Markov evolution corresponding to the two point motions $\{(X(t,t_0,\om, x), X(t,t_0,\om,y)): ( x,y)\in \R^{2d},\;  t\geqslant t_0\}$ by $P^{(2)}(t_0,(x,y); t, .)$ and $T^{(2)}_{t,t_0}$ respectively, and are defined by
\begin{align*}
P^{(2)}(t_0,(x,y); t, E) := \p\{\om: (X(t,t_0,\om,x), X(t,t_0,\om,y))\in E\}, \quad t\geqslant t_0, \quad E\in \Bb(\R^{2d})
\end{align*}
and
\begin{align*}
T^{(2)}_{t, t_0}h(x,y) &:= \int_{\R^{2d}}h(u,z)P^{(2)}(t_0,(x,y); t, dz\otimes du)\\ &=\E[h(X(t,t_0, \om,x), X(t,t_0,\om,y))],\quad h\in \mathcal{C}_b(\R^{2d}).
\end{align*}  
For $V \in \mathcal{C}_b^{1,2}(\R\times\R^{2d}),$ the generator $\mathcal{L}^{(2)}$ is defined by
\begin{align*}
\mathcal{L}^{(2)}V(t_0,x,y) = \lim_{t\downarrow t_0}\frac{\E[V(t, X(t, t_0, ., x),X(t, t_0, .,y))] - V(t_0, x,y)}{\vert t-t_0\vert}
\end{align*}
and, by It\^o's formula, we have $\mathcal{L}^{(2)}$ in the form
\begin{align*}
\mathcal{L}^{(2)}:= \frac{\partial  }{\partial t}+ \sum_{i=1}^{d}f_0^i(t,x)\frac{\partial}{\partial x_i}+\sum_{i=1}^df_0^i(t,y)\frac{\partial}{\partial y_i}+\frac{1}{2}\sum_{i,j=1}^{2d} \sum_{k=1}^mf_k^i(t,x,y)f_k^j(t,x,y)\frac{\partial^2}{\partial x_i\partial y_j},
\end{align*}
where $f_k^{i}(t,x,y) := f_k^{i}(t,x)$ and $f_k^{i+d}(t,x,y) := f_k^j(t,y), $ $i=1,2,\cdots,d$.

In particular, considering the difference between two solutions starting from two different initial values $x, y\in \Rd$ i.e., the process $\{X(t,t_0,\om,x)-X(t,t_0,\om,y):  x\neq y, t\geqslant t_0\}$,  the two-point generator $\mathcal{L}^{(2)}$ simplifies to 
\begin{align}\label{2pp}
\mathcal{L}^{(2)}V(t,x-y) =V_t(t,x-y)&+ V_x(t,x-y)\bigg(f_0(t,x)-f_0(t,y)\bigg)\nonumber\\ &\hspace{-1cm}+\frac{1}{2}\text{trace}\bigg( [\sigma(t,x)-\sigma(t,y)]^T\mathcal{H}V(t,x-y)[\sigma(t,x)-\sigma(t,y)]\bigg),
\end{align}
where $\sigma(t,x)= (f_1(t,x),f_2(t,x),\cdots, f_m(t,x)),$ $V_x = \Big(\frac{\partial V}{\partial x_i}\Big)_{ 1\leqslant i\leqslant d}$ and $ \mathcal{H}V = \Big(\frac{\partial^2 V}{\partial x_i \partial x_j}\Big)_{1\leqslant i,j\leqslant d}.$

Some stability and convergence results for stochastic flows from nontrivial reference solutions on a smooth manifold are normally via two-point generator (e.g., \cite{Bax2, Bax3, Har1}). Evidently, the two point generator $\mathcal{L}^{(2)}$ of the process $\{X(t,t_0,\om,x)-X(t,t_0,\om,y): x\neq y, t\geqslant t_0\}$ as $y\rightarrow x,$ is related to the one-point generator $\delta\mathcal{L}$ of the derivative flow $\{DX(t,t_0,\om, x)v: t\geqslant t_0\},$ (e.g., \cite{Uarnold, Arnold, Bax2, Hans3}). The consideration of two-point generator $\mathcal{L}^{(2)}$ is also one of the ways of making SDEs amenable to smooth ergodic theory. Schmalfuss in \cite{Schmalfuss} studied the existence of random stationary solutions and random attractors by Lyapunov function defined via two-point generator $\mathcal{L}^{(2)},$ our argument in this section, follows similar ideas in \cite{Schmalfuss}.

 Before presenting our results, we recall a variant of Doob's local martingale inequality that will be useful in the proofs.
\begin{prop}[Exponential local martingale inequality (cf. \cite{Mao})]\;\label{expl}\rm
Let $M= (M_t)_{t\ge 0}$ be a continuous local martingale. Then for any positive constants $\tau, \gamma, \delta,$ we have 
$$\p\bigg\{\om: \sup_{t\leq \tau}\big(M_t-\frac{\gamma}{2}\big<M\big>_t\big)>\delta\bigg\}\leq \exp(-\gamma\delta).$$
  
In particular, let $(M_t)_{t\ge 0}$ be a continuous real-valued local martingale vanishing at $t=0,$ $(\tau_k)_{k\ge 1}$ and $(\gamma_k)_{k\ge 0}$ be two sequences of positive numbers with $\tau_k\rightarrow\infty,$ $g(t)$ be a positive increasing function on $\Rp$ such that $$\sum_{k=1}^\infty g(k)^{-\theta}<\infty, \quad \text{for some }\theta>1.$$ Then, for almost all $\om\in \Om$ there is a random integer $k_0(\om)$ such that for all $k\ge k_0(\om)$ $$M_t\leq \frac{\gamma_k}{2}\big<M\big>_t+\frac{\theta}{\gamma_k}\log(g(k))\quad \text{on}\quad 0\leq t\leq \tau_k.$$  
\end{prop}
\begin{theorem}[Existence of random periodic solution]\;\label{tdper}\rm 
Let $X(t,t_0,\om,.)$ be a stochastic flow of diffeomorphisms induced by the SDE (\ref{phnl}) and suppose that $f_0,f_1,f_2,\cdots,f_m$ are periodic in $t$ with period $\tau>0.$
Let $V\in C^{1,2}(\Rp\times\Rd; \Rp)$ with $V(t,0)=0.$ Suppose there exist a function $\lambda: \R\rightarrow\R$ such that 
\begin{equation}\label{ilam}
\limsup_{t\rightarrow \infty}\frac{1}{2t}\int_{t_0}^t\lambda(s)ds<\alpha<0,
\end{equation}
\begin{align}\label{TWP}
\mathcal{L}^{(2)}V(t,x-y)\leq \lambda(t) V(t,x-y), \quad \text{and} \quad \vert x\vert^p\leq V(t,x), 
\end{align}
for all $t\in\R,$ $x,y\in\Rd,$  $p\ge 1.$
Suppose further that
 \begin{align}\label{Tempered}
\E\sup_{t>t_0}\ln V(t_0, X(t,t_0,\om,x)-x) <\infty.
\end{align}
 Then, there exists an $\F_{-\infty}^t$-measurable random variable $S(t,\om)$ which is the random periodic solution in the sense of Definition \ref{Chu}.
\end{theorem} 
\vspace{-.4cm}
\noindent {\it Proof.} The idea is to show that $\{X(t,t-n\tau,\om,x)\}_{n\in\N}$ is a Cauchy sequence in a space of continuous function $C([t_0,\infty); \Rd).$ For this, we shall first modify the initial value by a random variable in such a way that $X(t,t_0,\om,x)\neq X(t,t_0,\om,y)$ for some $y\in \Rd$ and all $t>t_0.$ This modification is necessary to avoid difficulties in the definiteness of some integrals and logarithms. We can disregard this modification at the end, using the fact $X(t,t_0,\om,x)$ is a homeomorphism (e.g., \cite{Hkunita, Kunita}) so that $X(t,t_0,\om,x)=X(t,t_0,\om,y)$ if and only if $x=y,$ (see \cite{Schmalfuss} for the case of time homogeneous flows, where such modification was made explicitly).    

Denote $X^x(t) := X(t,t_0,\om,x),\;\; X^y(t):=X(t,t_0,\om,y),$ we apply It\^o's formula on $\ln(V(t,X^x(t)-X^y(t)))$ to get
\begin{align*}\ln V(t,X^x(t)-X^y(t)) = \ln V(t_0,x-y)+\int_{t_0}^{t}\frac{\mathcal{L}^{(2)}V(s,X^x(s)-X^y(s))}{V(s,X^x(s)-X^y(s))}ds +N(t,\om) -\frac{1}{2}q(t,\om),
\end{align*} 
\newcommand{\HC}{\mathcal{G}}
where \begin{align*}
N(t,,\om)&=\int_{t_0}^t\big(\HC V(s,X^x(s),X^y(s))\big)dW_s, \quad q(t,\om)=\int_{t_0}^t\big(\HC V(s,X^x(s),X^y(s))\big)^2ds,\\
\HC V(t,x,y)&= \frac{V_x(t,x-y)\big(g(t,x)-g(t,y)\big)}{V(t,x-y)}, \quad W_t= (W_t^1,\cdots, W_t^m).
\end{align*}
We observe that $q(t,\om)$ is a quadratic variation of $N(t,\om)$ and as a consequence of Proposition \ref{expl}, we have that $$\p\bigg\{\om\in \Om: \sup_{t\in[t_0,t_0+k]}\bigg(N(t,\om)-\frac{1}{2}q(t,\om)\bigg)>2\ln k\bigg\}\leq \frac{1}{k^2}.$$ 
An application of  Borel-Cantelli lemma yields that for almost all $\om\in\Om,$ there exists a random integer $n_0=n_0(\om)$ such that for any $n\ge n_0$ 

$$\sup_{t_0\leq t\leq t_0+n}\bigg(N(t,\om)-\frac{1}{2}q(t,\om)\bigg)\leq 2\ln n.$$
In particular, \begin{equation}
\frac{1}{t}\bigg(N(t,\om)-\frac{1}{2}q(t,\om)\bigg)\leq \frac{2\ln n}{t_0+n-1}, \quad \text{for  $t_0+n-1$}\leq t\leq t_0+n.
\end{equation}
Next, applying the assumptions on $V$ and on the two point generator (\ref{TWP}), then the integrability condition (\ref{ilam}), we have that
\begin{align*}& \frac{1}{k-1}\sup_{n-1\leq t\leq n}\ln\vert X^x(t)-X^y(t)\vert \\ &\hspace{.5cm}\leqslant \frac{1}{p(n-1)}\sup_{n-1\leq t\leq n}\ln V(t, X^x(t)-X^y(t))\\ & \hspace{.5cm}\leqslant \frac{1}{p(n-1)}\ln V(t_0,x-y)+\frac{1}{p(n-1)}\sup_{n-1\leq t\leq n}\int_{t_0}^t\frac{\mathcal{L}^2V(s,X^x(s)-X^y(s))}{V(s,X^x(s)-X^y(s))}ds\\&\quad \quad\quad +\frac{1}{p(n-1)}\sup_{n-1\leq t\leq n}\bigg( N(t,\om)-\frac{1}{2}q(t,\om)\bigg)\\
 &\hspace{.5cm}\leqslant\ \frac{1}{p(n-1)}\ln V(t_0, x-y)  + \sup_{t\in [n-1,n]}\frac{1}{p(n-1)}\int_{t_0}^t\lambda(s)ds+ 2\frac{\ln n}{p(n-1)},\\
& \hspace{.5cm}\leqslant \frac{1}{p(n-1)}\ln V(t_0, x-y) +\frac{\alpha(n-1)}{p(t_0+n)}+2\frac{\ln n}{p(t_0+n-1)}
 \end{align*}
So, for $k$ large enough and for $\varepsilon>0,$ we have that
\begin{equation}\label{tpas}
\vert X(t,t_0,\om,x)-X(t,t_0,\om,y)\vert\leq  \left(V(t_0, x-y)\right)\exp({(\frac{\alpha}{p}+\varepsilon)(t-t_0)}), \quad \text{for all $t$}\ge t_0.
\end{equation}

Let $n\ge N,$  for $0<\varepsilon<-\frac{\alpha}{2p}$ and by flow property we have that
\begin{align*}&\frac{1}{n\tau}\ln \vert X(t,t-n\tau,\om,x)- X(t,t-N\tau,\om,x)\vert\\
&\hspace{2cm} =\frac{1}{n\tau}\ln \vert X(t,t-n\tau,\om,x)- X(t,t-n\tau,\om,X(t-n\tau, t-N\tau,\om,x)\vert\\
&\hspace{3cm}\leqslant \frac{1}{pn\tau}\ln V(t-n\tau, x- X(t-n\tau,t-N\tau,\om,x))+\frac{\alpha}{p}+\varepsilon.
\end{align*}
From the assumption that $\E\sup_{t>t_0}\ln V(t_0,X(t,t_0,\om,x)-x)<\infty,$ we have that for almost all $\om\in\Om,$ there exist a random variable $\beta(\om)$  with $\E[\vert \beta(\om)\vert]<\infty$ such that 
\begin{equation}\label{cacs}
\sup_{t\ge t_0}\vert X(t,t-n\tau,\om,x)-X(t,t-N\tau,\om,x)\vert \leqslant \beta \exp(\frac{\alpha n\tau}{2p}), \;\; \textrm{a.s.}.
\end{equation}
 So, the sequence $\{X(t,t-n\tau,\om,x)\}_{n\in\N}$ is Cauchy in the space $C([t_0, \infty);\R^d).$ Let  $S(t,\om)$ be the limit of the sequence $\{X(t,t-n\tau,\om,x)\}_{n\in\N},$ then by the continuity of the stochastic flow map $(t, t_0, x)\mapsto X(t,t_0,\om,x),$ we deduce for  $\tau>0$ that 
\begin{align}\label{R_pp}
\notag X(t+\tau, t,\om, S(t,\om))&=X(t+\tau,t,\om,\lim_{n\rightarrow\infty}X(t,t-n\tau,\om,x))\\  &= \lim_{n\rightarrow\infty}X(t+\tau,t+\tau-(n-1)\tau,\om,x)=S(t+\tau,\om).
\end{align}
Now, since the cofficients  $f_k: \; 0\leqslant k\leqslant m,$ of our SDE are time periodic with period $\tau$, we have that
 \begin{align*}X(t,t-n\tau,\theta_{\tau}\om,x)= &x+\int_{t-n\tau}^{t}f_0(r,X(r,r-n\tau,\theta_{\tau}\om,x))dr\\ & \quad \quad + \int_{t-n\tau}^{t}\sum_{k=1}^mf_k(r,X(r,r-n\tau,\theta_{\tau}\om,x))dW^k_{r+\tau}(\om)\\=&x+\int_{t+\tau-n\tau}^{t+\tau}f_0(r,X(r-\tau,r-\tau-n\tau,\theta_{\tau}\om,x)dr\\&\quad \quad +\int_{t+\tau- n\tau}^{t+\tau}\sum_{k=1}^{m}f_k(r,X(r-\tau,r-\tau-n\tau,\theta_{\tau}\om,x))dW_r^k(\om), \end{align*} and \begin{align*} X(t+\tau,t+\tau-n\tau,\om,x) & =x+\int_{t+\tau-n\tau}^{t+\tau}f_0(r,X(r,r-n\tau,\om,x))dr\\& \quad \quad +\int_{t+\tau-n\tau}^{t+\tau}\sum_{k=1}^{m}f_k(r,X(r,r-n\tau,\om,x))dW_r^k(\om).\end{align*} 
 By uniqueness of solution of SDE and the invariance of $\theta$ under $\p,$ we have that $$X(t,t-n\tau,\theta_{\tau}\om,x) = X(t+\tau,t+\tau-n\tau,\om,x),\;\; \p-\textrm{a.s.},$$ and then taking limit of both sides, we have from equation (\ref{R_pp}) that 
 \begin{equation*}
 S(t,\theta_{\tau}\om)= S(t+\tau,\om), \;\; \p-\textrm{a.s.}
 \end{equation*}

Moreover, from the inequality (\ref{tpas}), we see that for any $\F^r_{-\infty}$-measurable random variable $\xi(r,\om)$  
\begin{equation}\label{ExppER}
\lim_{n\rightarrow\infty}\sup_{t\in[n-1,n]}\vert S(t,\om)-X(t,r,\om,\xi(r,\om))\vert =0,  \;\; \p-\textrm{a.s.}
\end{equation}  
 exponentially fast.
\qed

\bigskip

We have the following corollory that captures cases of nonautonomous SDEs in many applications. They also arise in many theoretical study in stochastic analysis and nonlinear dynamical systems.
\begin{cor}\label{perie12}\rm
 Let $\{X(t,t_0,\om,.): t\geqslant t_0\}$ be a stochastic flow of diffeomorphisms induced by the SDE (\ref{phnl}) with $\tau$-periodic coefficients in time such that 
\begin{equation}
\begin{cases} \big \langle f_0(t,x)-f_0(t,y), x-y\big\rangle\leqslant \beta(t)\vert x-y\vert^2, \\ \vert f_k(t,x)-f_k(t,y)\vert \leqslant L(t)\vert x-y\vert,\;\; 1\leqslant k\leqslant m,
\end{cases}
\end{equation}
 for all $x,y\in\Rd,\; t\in\R$ and for some integrable functions $\beta: \R\rightarrow\R$ and $L:\R\rightarrow\Rp,$ such that 
 \begin{equation}\label{3.28}
\limsup_{t\rightarrow\infty}\frac{1}{2t}\int_{t_0}^t\beta (s)ds<0,\ \ 
\limsup_{t\rightarrow\infty}\frac{1}{2t}\int_{t_0}^tL (s)ds<\infty.
\end{equation}
Then the stochastic flow $\{X(t,t_0,\om,.): t\geqslant t_0\}$ generated by (\ref{phnl}) has a random periodic solution. 
\end{cor}
{\it Proof.} From (\ref{3.28}), it is easy to know that there exists a real number $p>1$ such that function  
$\lambda(s):=\beta(s)+\frac{(p-1)}{2}mL^2(s)$ satisfies 
\begin{equation}\label{3.28a}
\limsup_{t\rightarrow\infty}\frac{1}{2t}\int_{t_0}^t\lambda (s)ds<0.
\end{equation}
Take $V(t,x) = \vert x\vert^p$ for some $p>1$ and compute $\mathcal{L}^2V(t,x-y),$
$$\frac{\partial V(t,x)}{\partial x^i} = px^i\big(\sum_{n=1}^d (x^n)^2\big)^{\frac{p}{2}-1} = px^i\vert x\vert^{p-2},\;\; \frac{\partial^2V(t,x)}{\partial x^{i}\partial x^{j}} = 2(\frac{p}{2}-1)x^{i}x^{j}\vert x\vert ^{p-4} +\delta_{i,j}p\vert x\vert^{p-2},$$ where $\delta_{i,j}$ is the Kronnecker symbol. We now have for $1\leqslant k \leqslant m,$  \begin{align*}
& \sum_{k=1}^{m}\sum_{i,j=1}^{d}(x^{i}-y^{i})(f_k^{i}(t,x)-f_k^{i}(t,y))(f_k^{j}(t,x)-f_k^{j}(t,y))(x^{j}-y^{j}) \leqslant  mL^2(t)\vert x-y\vert^4,\\ 
 &\sum_{k=1}^{m}\sum_{i=1}^{d}(f_k^i(t,x)-f_k^i(t,y))(f_k^i(t,x)-f_k^i(t,y))  \leqslant mL^2(t)\vert x-y\vert^2,
 \end{align*}    
 and                                                                      
\begin{align*}
 \mathcal{L}^{(2)}V(t,x-y)&= \sum_{i=1}^{d}\big(f_0^i(t,x)-f_0^i(t,y)\big)p(x^i-y^i)\vert x-y\vert^{p-2}\\ &\quad \quad + \frac{1}{2}\sum_{i,j=1}^{d}\bigg[\bigg( \sum_{k=1}^{m}(x^{i}-y^{i})(f_k^{i}(t,x)-f_k^{i}(t,y))(f_k^{j}(t,x)-f_k^{j}(t,y))(x^{j}-y^{j})\bigg)\times\\& \qquad\qquad  \big\{2(\frac{p}{2}-1)p\vert x-y\vert^{p-4}\big\} \\ & \qquad+ \bigg( \sum_{k=1}^{m}(f_k^{i}(t,x)-f_k^{i}(t,y))(f_k^{j}(t,x)-f_k^{j}(t,y))\delta_{i,j}p\vert x-y\vert^{p-2}\bigg)\bigg]\\ &\leqslant p\beta(t)\vert x-y\vert^p + p\frac{mL^2(t)}{2}(p-1)\vert x-y\vert^p= p\lambda(t)\vert x-y\vert^p.
\end{align*}
Take $\hat{\lambda}(t) = p\lambda(t),$ it then follows from condition (\ref{3.28}) that
\begin{align*}
\limsup_{t\rightarrow\infty}\frac{1}{2t}\int_{t_0}^t \hat{\lambda}(s)ds<0.
\end{align*}
As $V(t,x) = \vert x\vert^p$, a one-point motion argument leads to if $X_0\in L^p(\Om)$ then $X(t,t_0,.,X_0)\in L^p(\Om)$, it the follows that $\E[V(t_0,X(t,t_0,\om,x))]<\infty,$ which gives us the temperedness assumption (\ref{Tempered}) on the vandom variable $V(X(t,t_0,\om,x))$ (see Lemma 1 in \cite{Schmalfuss}). 
\qed

\begin{rem}\rm\;

(i) If the functions $\beta$ and $L$ are continuous and periodic with period $\tau$, then (\ref{3.28}) can be replaced by a simple condition
\begin{eqnarray}\label{zhao2018a}
\int _0^{\tau}\beta(s)ds<0.
\end{eqnarray}

(ii)  Random periodic processes arise naturally in stochastic dynamical models in climatology, neuroscience, economics, molecular dynamics, etc.~This is due to the nonlinearity of the underlying vector fields and the onset of time-dependent random invariant sets, even in the case of temporal homogeneous vector fields.
We proved the existence of stable random periodic solutions and ergodicity of periodic measures for dissipative stochastic stochastic differential equations. 
 The assumption (\ref{zhao2018a}) 
we imposed is given in the sense of average, which is weaker than pointwise dissipativity. This is
natural in some physical models and can be verified in terms of the coefficients. For 
 example the following stochastic differential equations
\begin{eqnarray}\label{zhao2018b}
dx=[(-1+\gamma\sin  t)x-\delta x^3]dt+dW(t),
\end{eqnarray}
where $W(t)$ is a one-dimensional Brownian motion and $\delta \geq 0$ and $\gamma\in \R $ are real numbers, 
has a random periodic solution of period $2\pi$ according to Corollary \ref{perie12}. This natural result has not been discovered before.

 \end{rem}

\section{ Ergodicity in the random periodic regime}\label{Sec4}
In this section, we discuss ergodicity in the random periodic regime by considering probability measures induced by random periodic solutions. As we noted in the introduction, RDS theory is a systematic mix of stochastic analysis and dynamical systems theory. From stochastic analysis perspective, invariant probability measures are investigated via Markov transition probability function. In this sense,  ergodic theory is based on the dynamics of Markov evolution. From dynamical systems point of view, one studies random invariant probability measures whose conditional expectation with respect to a subalgebra of $\F$ has one to one correspondence with the invariant measure of Markov evolution.

 Here, we are interested in capturing ergodicity of the transition probability function in the random perodic regime (ergodicity of the law of random periodic solutions).  PS-ergodicity \footnote{Ergodicity on Poincar\'e sections} is a new form of ergodicity for stochastic dynamical systems in the random periodic setting, recently developed by Feng and Zhao in \cite{Zhao2}. It gives a new perspective and generalised form of Poincar\'e-Bendixson theorem for stochastic dynamical systems. We would like to argue using the information provided to us by the two-point generator, that the law of random periodic solutions form a family of PS-ergodic periodic measures. 
 \begin{definition}[Periodic measure \cite{Zhao2}]\;\label{Rand_pm}\rm
Let $\mathbb{M}$ be a Polish space, a  measure $\mu: \R\rightarrow \mathcal{P}(\mathbb{M})$ is called a {\it periodic measure} of period $\tau$ on the phase space $(\mathbb{M},\Bb(\mathbb{M}))$ for the Markovian stochastic flow $\{X(t,s,\om, .): \;t\geqslant s\}$ if  for $B\in \Bb(\mathbb{M})$ we have that 
 \begin{equation}\label{10.11}
 \mu_{s+\tau} = \mu_s\\ \quad \text{and} \quad \mu_{t+s}(B) = \int_{\mathbb{M}}P(s,x; t+s, B)\mu_s(dx) = (T^*_{t,s}\mu_s)(B), \quad s\in\R, \quad t\in \Rp.
 \end{equation}
 It is called a periodic measure with minimal period $\tau,$ if $\tau>0$ is the smallest number such that (\ref{10.11}) holds. 
 \end{definition} 
% \begin{rem}\rm
% A related definition of periodic measure was considered by Hasminksii in (\cite{Hasm}) in the sense that the transition probability $P(s,x; t,.)$ is periodic in time, i.e.,  for $k\in \Z$
% \begin{align}\label{12.h}
%\mu_s &= T^*_{k\tau,s}\mu_s = \mu_{s+k\tau}.
% \end{align}
% The second condition of (\ref{10.11}) is redudant in the Hasminksii's definition (\ref{12.h}), a counterexample which satisfies (\ref{12.h}) but not (\ref{10.11}) is considered in \cite{Zhao2}. We will see shortly, that the law of random periodic solutions satisfy  (\ref{10.11}).
% \end{rem}
 \bigskip
 
 Let $S: \R\times\Om\rightarrow \Rd$ be the random periodic solutions of the stochastic flow $\{X(t,s,\om,.): t\geqslant s\}$. We consider the probability measure
\begin{align}
  \mu_s(A) := \left(\p\circ S^{-1}(s,.)\right)(A) = \p(\{\om: S(s,\om)\in A\}), \;\; A\in \Bb(\Rd).
 \end{align} 
Then the measure $\mu_s$ is $\tau$-periodic as
 \begin{align}\label{PM3}
\nonumber \mu_{s+\tau}(A)= \p\{\om: S(s+\tau,\om)\in A\}
&= \p\{\om: S(s, \theta_{\tau}\om)\in A\}\\
& = \p\{\om: S(s,\om)\in A\}
 =\mu_s(A),
 \end{align}
 Moreover, as it was shown in \cite{Zhao2}, $\mu_s$ satistifies (\ref{10.11}). Thus, the law of random periodic solution satisfies Definition (\ref{Rand_pm}). 
% One also constructs time average measure $\bar \mu$ from the measure  $\mu_s$ \cite{Zhao2},
% \begin{align}\label{4.26}
%\bar \mu(A) &= \frac{1}{\tau}\int_0^{\tau}\mu_s(A)ds.
 %\end{align}
 
 \bn[Poincar\'e section for transition probability \cite{Zhao2}]\rm
 The collection of subsets $\{L_s: s\geqslant 0\}\subset \Bb(\Rd),$  are called the Poincar\'e sections of the transition probability function $P(s,x;t,.)$ if
  \begin{align*}
 L_{s+\tau} = L_s,
 \end{align*}
 and for any $s\in [0,\tau), \; t\geqslant 0,$
 \begin{align*}
 P(s,x; s+t,L_{s+t}) =1, \quad x\in L_s.
 \end{align*}
 
 \en
 \begin{rem}\rm\; \label{minP}
The choice of Poincar\'e section is not unique, example $L_s= \Rd$ and $L_s= \text{supp}(\mu_s)$ satisfy the definition of Poincar\'e section. However, the family $L_s= \text{supp}(\mu_s): s\in\R\}$ is a minimal Poincar\'e section \cite{Zhao2} or $n_s\tau$-irreducible Poincar\'e section. To see this, fix $s\in [0,\tau)$ and any open set $A_s\subset \textrm{supp}(\mu_s)$ with $\mu_s(L_s\backslash A_s)>0$, we have for all $x \in L_s,$
\begin{align*}
P(s, x;s+ n_s\tau, A_s)<1, \quad n\in \N.
\end{align*}
This implies that $A_s$ is not a Poincar\'e section for the transition probability $P(s,x; s+n\tau,.), \; n\in \N$.
\end{rem}
 \bn[PS-ergodicity \cite{Zhao2}]\rm
 The family of $\tau$-periodic measures $(\mu_s)_{s\in\R}\subset \mathcal{P}(\Rd)$ is {\it PS-ergodic,} if for each $s\in [0,\tau),$ $\mu_s$ as an invariant measure of the Markov evolution $(T_{s+k\tau,s})_{k\in\N}$ at the integral multiples of the period on the Poincar\'e section $L_s$ is ergodic.
 \en
 
 In fact, a periodic measure $\mu: \R\rightarrow \mathcal{P}(\Rd)$ is PS-ergodic, if  for any $A\in \Bb(\Rd)$ with $A\subset L_s$ and $\T\subset [0, \tau),$ we have
\begin{align}\label{1Ps}
\lim_{\N\rightarrow \infty}\frac{1}{\tau}\int_{0}^{\tau}\int_{\Rd}\bigg\vert\int_{\T}\Big(\frac{1}{N}\sum_{n=0}^{N-1}P(s,x; t+n\tau, A)-\mu_t(A)\Big)dt\bigg\vert\mu_s(dy)ds = 0.
\end{align}
The equation (\ref{1Ps}) is the Krylov--Bogolyubov scheme for periodic measure \cite{Hasm, Zhao2}.
\medskip

Given the above preparation, we are now ready to prove the PS-ergodicity of periodic measures generated by a class of SDEs with time periodic coefficients satisfying conditions of Theorem \ref{tdper}.
Precisely, we want to prove the convergence of  Kyrlov--Bogolyubov scheme for periodic measures.  However, strong Feller property of the Markov evolution $(T_{s+n\tau,s})_{n\in \N}$ is crucial in the proof of the convergence of Krylov--Bogolyubov scheme. Recall that a Markov evolution  $(T_{s+n\tau,s})_{n\in\N}$ has strong Feller property if $\Bb_b(\Rd)\ni \varphi\mapsto T_{s+n\tau,s}\varphi \in \mathcal{C}_b(\Rd)$.
Equivalently, the Markov evolution $(T_{s+n\tau,s})_{n\in \N}$ has strong Feller property if and only if
\btm
\item[(a)] $(T_{s+n\tau,s})_{n\in \N}$ is a Feller semigroup, i.e., $T_{s, s+n\tau}: \mathcal{C}_b(\Rd)\rightarrow\mathcal{C}_b(\Rd)$ and 
\item[(b)] the family $\{P(s,x_m;s+n\tau,.): m\in \N\}$ is equicontinuous.
\etm  
The first item, follows from the existence of stochastic flow (see, e.g.,~\cite{Kunita, Hasm}), the second is not simple and would require further non-degeneracy condition on the diffusion coefficients of the SDE. To this end, we start with following  lemma.
%\vspace{-0.1cm}\bn[Pseudo-inverse of a matrix-valued vector  field]\rm\label{Pseudo}
\mbox{}
%\btm[leftmargin = 0.7cm]
%The psuedo-inverse $\tilde{\sigma}^{-1}$ of a matrix vector field $\sigma\,{:}\; \Rd\rightarrow L(\R^m;\Rd)$ is defined via the  limit  
%\begin{align*}
%\tilde{\sigma}^{-1} &= \lim_{\varepsilon \rightarrow 0}\left( \sigma^T\sigma+ \varepsilon\textrm{I}\right)^{-1}\sigma^T = \lim_{\varepsilon \rightarrow 0}\sigma^T\left( \sigma\sigma^T+ \varepsilon\text{I}\right)^{-1}.
%\end{align*}
%In a similar fashion, the pseudo-inverse of $a = \sigma\sigma^T$ is given by  $\tilde{a}^{-1} = \tilde{(\sigma^T)}^{-1}\tilde{\sigma}^{-1}$ %with the following properties
%such that
%\begin{equation*}
%\sigma = a\tilde{\sigma}^{-1} \;\; \textrm{and}\;\; \tilde{\sigma}^{-1} = \sigma^T \tilde{a}^{-1}.
%\end{equation*}

%Geometrically, for $a: \Rd\rightarrow L(\Rd; \Rd)$ and $\tilde{a}^{-1}: \Rd\rightarrow L(\Rd;\Rd)$, the vector spaces $\Rd$ and $L(\Rd;\Rd)$ can be decomposed as  $\Rd = [\textrm{ker} (a)]^{\bot}\oplus \textrm{ker}(a)$ and $L(\Rd;\Rd) = \textrm{ran}(a)\oplus[ \textrm{ran}(a)]^{\bot}.$ The restriction $\tilde{a}: [\textrm{ker}(a)]^{\bot}\rightarrow \textrm{ran}(a)$ is an isomorphism; therefore, the pseudo-inverse $\tilde{a}^{-1}$ is defined on $\textrm{ran}(a)$ to be the inverse of $\tilde{a},$ namely
%\begin{align*}
%\tilde{a}^{-1} &=
%\begin{cases}
%\;\;\;0^T, & \textrm{if}\;\;  a = 0,\\
%a^T(a^Ta)^{-1}, &\textrm{if}\;\; a\neq 0.
%\end{cases}
%\end{align*}
%\en

\begin{lem}\label{Strong Feller}\rm
Suppose there exists $\tau>0$ such that $f_k(t+\tau, x) = f_k(t,x), \; 0\leqslant k\leqslant m$ with $f_0\in \mathcal{C}_b^1(\R\times\Rd; \Rd)$ and $f_k\in \mathcal{C}_b^2(\R\times\Rd; \Rd), \; 1\leqslant k\leqslant m$ . If in addition to the assumptions of Theorem \ref{tdper}, there exists $C\geqslant 1$ such that $V(t,x)\leqslant C\vert x\vert^p,\; p\geqslant 1$ and there exists $K>0$ such that 
\begin{align}\label{Nondege}
\sup_{t\in [s, s+n\tau]}\Vert \sigma^{-1}(t,x)\Vert_{\text{HS}}\leqslant K, \qquad \forall x\in \Rd, \quad n\geqslant 1,
\end{align} 
where $\sigma^{-1}(t,x)$ is the right inverse of  $\sigma(t,x) := (f_1(t,x), \cdots, f_m(t,x))$.

Then, there exists $0<K_C<\infty$ such that for $x, y\in L_s$ with  $s\in [0, \tau),$ and $h\in \mathcal{C}_b(\Rd)$, we have
\begin{align}
\vert T_{s+n\tau,s}h(x) -T_{s+n\tau,s}h(y)\vert \leqslant \frac{K_C}{\sqrt{n\tau}}\exp\left(\frac{1}{2}\int_{s}^{s+n\tau}\lambda(r)dr\right)\Vert h\Vert_{\infty} \vert x-y\vert.
\end{align}
\end{lem}
\noindent {\it Proof.}
First, we show for any $Y ,Z\in L_p(\Om, \F_{-\infty}^s, \p)$ for $p\geqslant 1$ and $t\geqslant s$,
\begin{align}\label{Estim1}
\E\vert X(t, s,\om, Y(s,\om))- X(t, s,  \om, Z(s,\om))\vert^p\leq C\exp\Big(\int_s^{t}\lambda(r)dr\Big)\E\vert Y(s,\om)-Z(s,\om)\vert^p.
\end{align}
 For this, set $\alpha(t) = \exp\Big(-\int_s^t\lambda(r)dr\Big)$ and $M(t,s,\om,x) = \sum_{k=1}^m\int_{s}^tf_k(r, x)dW^k_r,$ then by It\^o's formula (Theorem 8.1 in \cite{Hkunita}) we have 
\begin{align*}
&d\bigg(\alpha(t)V(t, X(t, s, \om,Y)-X(t, s, \om, Z))\bigg) \\ &\hspace{1cm}=-\lambda(t)\alpha(t) V(t, X(t,s, \om,Y)-X(t, s, \om, Z))dt\\ & \hspace{2cm}
+\alpha(t)\mathcal{L}^{(2)}V(t, X(t,\om,Y)-X(t, s, \om, Z))dt\\
&\hspace{2cm} +\alpha(t)V(t, X(t,s,\om,Y)-X(t, s,\om,Z))d\Big(M(X(t,s, \om,Y)-M(X(t, s, \om,Z))\Big)\\
&\hspace{1cm}\leqslant \alpha(t)V(t, X(t,s, \om,Y)-X(t,s,\om,Z))d\Big(M(X(t, s, \om,Y)-M(X(t,s,\om, Z))\Big).
\end{align*}
This implies that for $t\geqslant s$,
\begin{align*}
\E\bigg(\vert X(t,s,\om,Y(s,\om))- X(t,s,\om,Z(s,\om))\vert^p \bigg)&\leqslant \E\bigg(V(t, X(t,s,\om,Y(s,\om)) -X(t,s,\om,Y(s,\om))\bigg)\\
&\hspace{.5cm}\leqslant C\exp\Big(\int_s^{t}\lambda(r)dr\Big)\E\vert Y(s,\om)-Z(s,\om)\vert^p.
\end{align*}
In particular, for $x\neq y\in L_s, \;s\in [0,\tau),$ we obtain
\begin{align}
\E\vert X(s+n\tau,s, \om, x) - X(s+n\tau, s, \om, y)\vert^p\leq C\exp\Big(\int_s^{s+n\tau}\lambda(r)dr\Big)\E\vert x-y\vert^p.
\end{align}
Since the coefficients $(f_k)_{k=0}^m$ are such that the derivative flow $v_{s+n\tau,s}:=D_xX(s+n\tau,s,\om, x)v$ at $x\in L_s$ in the direction $v\in \Rd$ exists for almost all  $\om\in\Om$ and $s\in [0, \tau),\; n\in \N,$ then as $y\rightarrow x,$ we have
\begin{align*}
\E\vert v_{s+n\tau,s}\vert=\E\vert D_xX(s+n\tau,s,\om,x)v\vert\leqslant C^{1/p}\exp\left(\int_s^{s+n\tau}\lambda(r)dr\right)^{1/p}\vert v\vert.
\end{align*}
This implies that 
\begin{align}\label{Df_bound1}
\E\Vert D_xX(s+n\tau,s,\om,x)\Vert = \E\left(\sup_{v\neq 0}\frac{\vert v_{s+n\tau,s}\vert}{\vert v\vert}\right)\leqslant C^{1/p}\exp\left(\int_s^{s+n\tau}\lambda(r)dr\right)^{1/p}.
\end{align}
 Next, by It\^o's formula, we have for $h\in \mathcal{C}_b^2(\Rd),$
 \begin{align*}
 h(X(s+n\tau,s,\om,x)) = T_{s+n\tau,s}h(x) +\int_s^{s+n\tau}D_x(T_{s+n\tau,r}h)(X(r,s,\om,x))\sigma(r, X(r,s,\om,x))dW_r.
 \end{align*}
  Multiplying both sides by $\int_s^{s+n\tau}\sigma^{-1}(r,X(r,s,\om,x))DX(r,s,\om,x)vdW_r,$ taking expectation and applying It\^o isometry and Fubini's theorem, we have
 \begin{align*}
&\hspace{-2cm} \E\left(  h(X(s+n\tau,s,\om,x)) \int_s^{s+n\tau}\sigma^{-1}(r, X(r,s,\om,x))D_xX(r,s,\om,x)vdW_r\right) \\ &= \E\left( \int_s^{s+n\tau} D_x(T_{s+n\tau,r}h)(X(r,s,\om,x))\cdot D_xX(r,s,\om,x)vdr\right)\\
 &=\int_s^{s+n\tau} D_x\E[T_{s+n\tau,r}h(X(r,s,\om,x))]vdr.
 \end{align*}
Using Markov property of the stochastic flow $\{X(t,s,\om,.): t\geqslant s\}$, we have that $$ \E[T_{s+n\tau,r}h(X(r,s,\om,x))] = T_{r,s}\circ T_{ s+n\tau,r}h(x) = T_{s+n\tau,s}h(x),$$ so that, for any $h\in \mathcal{C}_b^2(\Rd),$ we arrive at 
\begin{align}\label{ElwoLI}
D_xT_{s+n\tau,s}h(x)v = \frac{1}{n\tau} \E\Big[ h(X(s+n\tau,s,\om,x)\int_s^{s+n\tau}\!\!\!\! \sigma^{-1}(r,X(r,s,\om,x))v_{r,s} dW_r\Big],
\end{align} 
 Next, since $\mathcal{C}_b^2(\Rd)$ is dense in $\mathcal{C}_b^1(\Rd),$ we obtain a version of {\it Bismut-Elworthy-Li formula} (e.g., \cite{Da Prato2, ElwF}) (\ref{ElwoLI}) for all $h\in \mathcal{C}_b^1(\Rd)$.

Since $\mathcal{C}_b^1(\Rd)$ is dense in $\mathcal{C}_b(\Rd),$ we have $(h_m)_{m\in \N}\subset \mathcal{C}^1_b(\Rd)$ such that $h_m\rightarrow h\in \mathcal{C}_b(\Rd)$
and 
\begin{align*}
\lim_{m\rightarrow\infty}T_{s+n\tau,s}h_m(x) &= T_{s+n\tau,s}h(x),\\
\lim_{m\rightarrow \infty}D_xT_{s+n\tau,s}h_m(x)\cdot v& = \frac{1}{n\tau} \E\Big[ h(X(s+n\tau,s,\om,x)\int_s^{s+n\tau}\!\!\sigma^{-1}(r, X(r,s,\om,x)))\cdot v_{r,s} dW_r\Big].
\end{align*}
convergence being uniform (cf.~\cite{Da Prato2}).  On the other hand, since $f_0\in \mathcal{C}_b^1(\R\times\Rd;\Rd), \; f_k\in \mathcal{C}_b^2(\R\times\Rd;\Rd), \; 1\leqslant k\leqslant m$ with the non-degeneracy condition (\ref{Nondege}), there exists a function $0<\rho_{s+n\tau,s}\in \mathcal{C}_b^1(\Rd)\times\mathcal{C}_b^1(\Rd)$ such that $P(s,x;s+n\tau,dy) = \rho_{s+n\tau,s}(x,y)dy$ (e.g.,~\cite{Da Prato2, Hasm, Kunita}). This implies that 
\begin{align*}
\notag \lim_{m\rightarrow \infty}D_xT_{s+n\tau,s}h_m(x)\cdot v &= \lim_{m\rightarrow\infty}\int_{\Rd}h_m(y)D_x(\rho_{s+n\tau,s})(x,y)\cdot vdy\\
\notag &=\int_{\Rd}h(y)D_x\rho_{s+n\tau,s}(x,y)\cdot v dy=D_x(T_{s+n\tau,s}h)(x), \quad 
\end{align*} 
so (\ref{ElwoLI}) holds for all $h\in \mathcal{C}_b(\Rd).$

Next, by the equality (\ref{ElwoLI}), Cauchy--Schwartz inequality, It\^o isometry and the condition (\ref{Nondege}), we have 
\begin{align}\label{Df_bound2}
\notag |D_xT_{s+n\tau,s}h(x) v|^2&\leqslant  T_{s+n\tau,s}h^2(x)\frac{1}{(n\tau)^2}\E\left(\int_s^{s+n\tau}\vert \sigma^{-1}(r,X(r,s,\om,x))\cdot v_{r,s}\vert^2 dr\right)\\ 
&\leqslant \Vert h\Vert^2_{\infty} \frac{K^2}{n\tau}\E\Vert D_xX(s+n\tau,s,\om,x)\Vert^2\vert v\vert^2
\end{align}
Comparing (\ref{Df_bound2}) with (\ref{Df_bound1}), we have  
\begin{align}
\vert D_x T_{s+n\tau,s}h(x) v\vert\leqslant \frac{K_C}{\sqrt{n\tau}}\exp\left(\frac{1}{2}\int_s^{s+n\tau}\lambda(r)dr\right)\Vert h\Vert_{\infty}\vert v\vert, \qquad  x\in L_s,\; v \in \Rd, 
\end{align}
where $K_C= K\sqrt{C}$.
Finally, let $\eta^\ell(x,y) = \ell x+(1-\ell)y,\; x, y\in L_s, \; \ell\in [0, 1]$, for $h\in \mathcal{C}_b(\Rd)$, the mean value theorem, leads to 
\begin{align*}
\vert T_{s+n\tau,s}h(x)- T_{s+n\tau,s}h(y)\vert &=\vert \int_0^1D_x(T_{s+n\tau,s}h)(\eta^\ell(x,y))\cdot (x-y) d\ell\vert \\ &\hspace{1cm}\leqslant \frac{K_C}{\sqrt{n\tau}}\exp\left(\frac{1}{2}\int_{s}^{s+n\tau}\lambda(r)dr\right)||h||_{\infty}\vert x-y\vert.
\end{align*}
%It then follows that for  $h\in \mathcal{C}_b(\Rd),$
%\begin{align*}
%\vert T_{s+n\tau,s}h(x)- T_{s+n\tau,s}h(y)\vert & \leqslant \frac{K_{C}}{\sqrt{n\tau}}\exp\left(\frac{1}{2}\int_s^{s+n\tau}\lambda(r)dr\right)||h||_{\infty}\vert x-y\vert.
%\end{align*}
\qed

\begin{theorem}\label{Ps_ergodic}\rm
Suppose the conditions of Lemma \ref{Strong Feller} are satisfied.  Moreover, assume there exist $\tau$-periodic function $V\in \mathcal{C}^{1,2}(\R\times\Rd;\Rp)$ satisfying the conditions of Theorem \ref{tdper} and $C\geqslant 1$ such that $V(t,x)\leqslant C\vert x\vert^p$ and 
\begin{align}\label{Strenh}
\E\big[\sup_{t\geqslant t_0} V(t_0, X(t,t_0,\om,x) -x)\big]<\infty
\end{align}
Then, the family of periodic measures $(\mu_s)_{s\in [0, \tau)}$ induced by the random periodic path $S(s,\om)$ is PS-ergodic.
\end{theorem}
\noindent {\it Proof.} 
{\bf Step I:}
First, we show that there exists $0<\tilde{K}<\infty$ such that for any initial value $x\in L_p(\Om, \F_{-\infty}^s,\p), \; p\geqslant 1$, we have 
\begin{align}\label{Estima2}
\E\bigg(\vert X(s+n\tau, s, \om,x)-S(s+n\tau,\om)\vert^p\bigg)\leqslant \tilde{K}\exp\left(\int_{s}^{s+n\tau}\lambda(u)du\right), \quad \forall n\in \N.
\end{align}
For this, we recall from the definition of random periodic solution that $S(s+n\tau, \om) = X(s+n\tau,s,S(s,\om)), \; \p-\text{a.s.,}$ then the estimate (\ref{Estim1}) yields,
\begin{align*}
\E\bigg(\Big\vert X(s+n\tau, s, \om,x)-S(s+n\tau,\om)\Big\vert^p\bigg) &= \E\bigg(\Big\vert X(s+n\tau, s, \om,x)-X(s+n\tau, s, \om, S(s,\om))\Big\vert^p\bigg)\\
&\hspace{.5cm} \leqslant C\exp\left(\int_s^{s+n\tau}\lambda(u)du\right)\E\vert x-S(s,\om)\vert^p.
\end{align*} 
Next, recall from the construction in Theorem \ref{tdper} that $S(t,\om)$ is the limit of the Cauchy sequence $\{X(s, s-n\tau, \om,x): n\in \N\}$ in $\mathcal{C}([s,\infty);\Rd)$  and $\E\vert S(s,\om)-X(s,s-n\tau,\om,x)\vert$ converges to zero exponentially according to equation (\ref{ExppER}) or (\ref{cacs}). In view of these and the condition $\vert x\vert^p\leqslant V(t,x)$ with time periodicity of $V(t,x),$ together with triangle inequality, we have
\begin{align*}
&\hspace{-1.5cm}\E\bigg(\vert X(s+n\tau, s, \om,x)-S(s+n\tau,\om)\vert^p\bigg)\\ &\leqslant C\exp\left(\int_s^{s+n\tau}\lambda(u)du\right)\bigg\{\E\vert x-X(s, s-n\tau, \om,x)\vert^p + \beta\exp( \frac{\alpha n\tau}{2})\bigg\}\\
&\leqslant C\exp\left(\int_s^{s+n\tau}\lambda(u)du\right)\bigg\{\E[V(s,X(s, s-n\tau,\om,x)-x)]+ \beta \exp(\frac{\alpha n\tau}{2})\bigg\}\\
&= C\exp\left(\int_s^{s+n\tau}\lambda(u)du\right)\bigg\{\E[V(s-n\tau,X(s, s-n\tau,\om,x)-x)]+\beta\exp(\frac{\alpha n\tau}{2})\bigg\},
\end{align*}
where $\beta$ and $\alpha$ are constants from the inequality (\ref{cacs}).\\
Next, by the condition (\ref{Strenh})
we have a random variable $\gamma$ with $\E[\gamma(\om)]<\infty$ such that 
\begin{align*}
\E\bigg(\vert X(s+n\tau, s, \om,x)-S(s+n\tau,\om)\vert^p\bigg)&\leqslant  \tilde{C}\E[\gamma(\om)]\exp\left(\int_s^{s+n\tau}\lambda(u)du\right)\\
&\hspace{.5cm} = \tilde{K} \exp\left(\int_s^{s+n\tau}\lambda(u)du\right).
\end{align*}
{\bf Step II:}  In this step, we show that there exists  $K_1>0$ such that 
\begin{align}\label{Lipq1}
\Big\vert T_{s+n\tau,s}\varphi(x) - \int_{L_s}\varphi d\mu_s\Big\vert\leqslant K_1\Vert \varphi\Vert_{\infty} \exp\left(p^{-1}\int_s^{s+n\tau}\lambda(r)dr\right), \quad \forall\; \varphi \in \mathcal{C}_b(L_s).
\end{align}
To this end,  we start by recalling that the periodic measure $\mu_s$ is invariant under the Markov evolution $(T_{s+n\tau, s})_{k\in\N}$, so that for any $h\in \text{Lip}_b(L_s)$
\begin{align}\label{Lipq}
\left\vert T_{s+n\tau,s}h(x)-\int_{L_s}hd\mu_s\right\vert 
\notag &= \left\vert \int_{L_s}\Big(T_{s+n\tau,s}h(x)-T_{s+n\tau,s}h(y)\Big) \mu_s(dy)\right\vert\\
\notag &\leqslant \Vert h\Vert_{\text{BL}}\int_{L_s}\E\left\vert X(s+n\tau,s,\om,x) - X(s+n\tau,s,\om,y)\right\vert
 \mu_s(dy)\\
 &\leqslant \tilde{K} \Vert h\Vert_{\text{BL}}\exp\left( p^{-1}\int_s^{s+n\tau}\lambda(u)du\right),
\end{align}
where we have applied H\"older's inequality and (\ref{Estima2}) in the last line.
Let $\varphi\in \mathcal{C}_b(L_s)$ be given. Setting $h = T_{s+\tau,s}\varphi $ in (\ref{Lipq}), then, by Lemma \ref{Strong Feller}
and by the invariance of $\mu_s$ under the Markov evolution $(T_{s+k\tau,s})_{k\in \N},$ we have 
\begin{align*}
\left\vert T_{s+\tau+n\tau,s}\varphi(x)- \int_{L_s}T_{s+\tau,s}\varphi d\mu_s\right\vert = \left\vert \int_{L_s}\left(T_{s+\tau+ n\tau,s}\varphi(x) -T_{s+n\tau +\tau,s}\varphi(y)\right)\mu_s(dy)\right\vert\\
\leqslant \tilde{K}\Vert T_{s+\tau,s}\varphi\Vert_{\text{BL}}\exp\left(p^{-1}\int_s^{s+n\tau}\lambda(r)dr\right)\\
\leqslant K_{\tau}\Vert \varphi\Vert_{\infty}\exp\left(p^{-1}\int_s^{s+n\tau}\lambda(r)dr\right),
\end{align*}
where $K_\tau = \frac{\tilde{K} K_C}{\sqrt{\tau}}\exp\left(\frac{1}{2}\int_s^{s+\tau}\lambda(r)dr\right)$.
It then follows that for any $\varphi \in \mathcal{C}_b(L_s)$ and $n>1,$ 
\begin{align*}
\left\vert T_{s+n\tau,s}\varphi(x)- \int_{L_s}\varphi d\mu_s\right\vert\leqslant K_{\tau}\Vert \varphi\Vert_{\infty}\exp\left(p^{-1}\int_{s}^{s+n\tau-\tau}\lambda(r)dr\right),
\end{align*}
this implies that there exists $K_1>0$ such that (\ref{Lipq1}) holds.\\
%Finally, since 
%\begin{align*}
%\Big\vert T_{s+n\tau,s}\varphi(x) - \int_{L_s}\varphi d\mu_s\Big\vert \leqslant 2\Vert \varphi \Vert_{\infty}, \quad \forall\; \varphi\in \mathcal{C}_b(L_s),
%\end{align*}
%then,  there exist $K_1>0$ such that 
%\begin{align*}
%\Big\vert T_{s+n\tau,s}\varphi(x) - \int_{L_s}\varphi d\mu_s\Big\vert\leqslant K_1\Vert \varphi\Vert_{\infty} \exp\left(p^{-1}\int_s^{s+n\tau}\lambda(r)dr\right), \quad \forall \; \varphi \in \mathcal{C}_b(L_s).
%\end{align*}
{\bf Step III:}
To complete the proof, we employ density argument and use \ref{Lipq1} to obtain the convergence of Krylov--Bogolyubov scheme for periodic measure.
Now, let $A_s\subset L_s$ be a closed set,  take $\varphi = \I_{A_s}$ and consider the sequence of functions $(\varphi_m)_{m\in \N}$ defined by 
\begin{align*}
\varphi_m(x) = \begin{cases} 1, & \quad \text{if}\quad x\in A_s,\\
1-2^md(x, A_s), & \quad \text{if}\quad d(x,A_s)\leqslant2^{-m},\\
0, &\quad \text{if} \quad d(x,A_s)\geqslant 2^{-m},
\end{cases}
\end{align*}
where $d(x,A_s) = \inf\{\vert x-y\vert: y\in A_s \},\; x\in L_s$. Then,
\begin{align*}
\varphi_m(x)\rightarrow \varphi(x), \quad \text{as}\; m\rightarrow\infty, \quad x\in L_s.
\end{align*}
Now, for all $s\in [0,\tau)$ we have 
\begin{align*}
T_{s+n\tau,s}\varphi_m(x)\rightarrow T_{s+n\tau,s}\varphi(x) = T_{s+n\tau,s}\I_{A_s}(x), 
\end{align*}
this implies that $T_{s+n\tau,s}\varphi_m\in \mathcal{C}_b(L_s)$ and as $\mu_s$ is invariant under $T_{s+n\tau,s},$ then (\ref{Lipq1}) leads to 
\begin{align}\label{Df_bound3}
\notag \left\vert P(s,x;s+n\tau,A_s) -\mu_s(A_s)\right\vert&=\left\vert T_{s+n\tau,s}\I_{A_s}(x)-\mu_s(A_s)\right\vert\\ &\hspace{1cm} \leqslant K_1 \exp\left(p^{-1}\int_s^{s+n\tau}\lambda(u)du\right).
\end{align}
By covering lemma (e.g.,~\cite{AMB00}), the inequality (\ref{Df_bound3}) holds for any $A\in \Bb(\Rd)$  with $A\subset L_s,$  thus, for $\T\subset [0, \tau),$ we have  
\begin{align*}
\int_{\T} \left\vert P(s,x; s+n\tau, A)-\mu_s(A)\right\vert ds&\leqslant \int_{0}^{\tau} \left\vert P(s, x;s+n\tau,A) -\mu_s(A)\right\vert ds\\ &\leqslant K_1 \int_{0}^{\tau}\exp\left(p^{-1}\int_s^{s+n\tau}\lambda(u)du\right)ds\\
&= K_1\int_{0}^{\tau}\exp\left(\frac{1}{2pn\tau}\int_{s}^{s+n\tau}\lambda(u)du\right)^{2n\tau}ds.
\end{align*}  
Next, we use the Chapmann--Kolmogorov equation for tranisition probability to obtain
\begin{align*}
\bigg\vert \int_{\T}\left[P(s,x;t+n\tau,A)-\mu_t(A)\right]dt\bigg\vert &= \bigg\vert \int_{\T}\left[\int_{\Rd}P(t,y;t+n\tau,A)-\mu_t(A)\right]P(s,x;t,dy)dt\bigg\vert\\
&\leqslant\int_{0}^{\tau}\int_{\Rd} K_1\exp\left(\frac{1}{2pn\tau}\int_t^{t+n\tau}\lambda(u)du\right)^{2n\tau}P(s,x;t,dy)dt\\
&=\int_{0}^{\tau}K_1\exp\left(\frac{1}{2pn\tau}\int_t^{t+n\tau}\lambda(u)du\right)^{2n\tau}\int_{\Rd}P(s,x;t,dy)dt\\
&= K_1\int_0^{\tau}\exp\left(\frac{1}{2pn\tau}\int_t^{t+n\tau}\lambda(u)du\right)^{2n\tau}dt.
\end{align*} 
By the condition (\ref{ilam}) of Theorem \ref{tdper},  we have there exist $0<\beta<1,\; 0<K_2<\infty,$ such that 
\begin{align*}
\bigg\vert\int_{\T}\left(P(s,x;t+n\tau,A)-\mu_t(A)\right)dt\bigg\vert \leqslant \int_{\T}\left\vert P(s,x; t+n\tau, A)-\mu_t(A)\right\vert dt\leqslant K_2\beta^{n\tau}.
\end{align*}
It then follows that
 \begin{align*}
\frac{1}{\tau}\int_0^{\tau}\int_{\Rd}\left\vert \int_{\T}\left\{\frac{1}{N}\sum_{n=0}^{N-1}P(s,x;t+n\tau, A) - \mu_t(A)\right\}dt\right\vert\mu_{s}(dx)ds\leqslant K_2\sum_{k=0}^{N-1}\beta^{n\tau}.
\end{align*}
This implies the convergence of the Krylov-Bogolyubov scheme for periodic measures.
\qed

 \subsection*{Acknowledgements} HZ would like to acknowledge the financial supports of an EPSRC Established Career Fellowship (ref: EP/ 50053293/1) and the Royal Society Newton Fund (ref: NA150344).


\begin{thebibliography}{plain}
\setlength{\itemsep}{-1mm}
\bibitem{AMB00} L. Ambrosio, N. Fuso and D. Pallara, \emph{Functions of bounded variation and free discontinuity problems}, Clarendon Press, Oxford, 2000.
\bibitem{Uarnold} L. Arnold, \emph{The unfolding dynamics in stochastic analysis,} Computational and Applied Mathematics, Vol.16 (1997), 3 - 25.
\bibitem{Arnold} L. Arnold, \emph{Random Dynamical Systems,} Springer, 1998.
\bibitem{Arnoldp} L. Arnold and M. Scheutzow, \emph{Perfect cocycles through stochastic differential equations,} Probability Theory and Related Fields, Vol. 101 (1995), 65 - 88.
%\bibitem{Bax04} P. Baxendale, \emph{Stochastic averaging and asymptotic behaviour of the stochastic Duffing--van der Pol equation}, Stochastic processes and their applications, 113 (2004) 235 -- 272.
\bibitem{Bax2} P. Baxendale, \emph{Stability and Equilibrium properties of stochastic flows of diffeomorphisms,} In: Diffusion processes and related problems in analysis, Birkh\"auser, 1992.
\bibitem{Bax3} P. Baxendale, \emph{Statistical equilibrium and two-point motion for a stochastic flow of diffeomorphisms,} In: Spatial stochastic processes, Birkh\"auser, 1991.
%\bibitem{Bax4} P. Baxendale, \emph{The Lyapunov spectrum of a stochastic flow of diffeomorphisms}, In: Lyapunov exponents, Springer, 1984.
\bibitem{Bax} P. Baxendale, \emph{Wiener processes on manifolds of maps,} Proceedings of the Royal Society of Edinburgh, Section A, Vol. 87 (1980), 127 - 152.
\bibitem{Bism} J. M. Bismut, \emph{A generalized formula of Ito and some other properties os stochastic flows,} Z. Wahrscheinlinlichkeeitstheorie verw. Gebiete, Vol. 55 (1981), 331 - 350.
%\bibitem{Rockner} V. I. Bogachev, N. V. Krylov, M. R\"ockner and S. V. Shaposhnikov, \emph{Fokker-Planck-Kolmogorov Equations}, American Mathematical Society, 2015.
\bibitem{Caraballo} T. Caraballo and P. E. Kloeden and B. Schmalfuss, \emph{Exponentially Stable Stationary Solutions of Stochastic Evolution Equations and Their Perturbations,} Appl. Math. Optim. Vol. 50 (2004), 183 - 207.
\bibitem{Elw2} A. P. Carverhill and K. D. Elworthy, \emph{Flows of Stochastic Dynamical Systems: The Functional Analytic Approach,} Z. Wahrscheinlichkeitstheorie verw. Gebiete, Vol. 65 (1983), 245 - 267. 
%\bibitem{Floquet} C. Chicone, \emph{Ordinary Differential Equations with Applications,} Springer, 2006.
\bibitem{Hans3} H. Crauel, \emph{Extremal Exponents of Random Dynamical Systems Do Not Vanish,} Journal of Dynamics and Differential Equations, Vol. 2 (1990), 245 - 291.
%\bibitem{Crau2} H. Crauel, \emph{Random Probability Measures on Polish spaces}, Taylor and Francis,  2002.
%\bibitem{Da Prato} G. Da Prato and J. Zabczyk, \emph{Ergodicity for infinite dimensional systems,} Cambridge University Press, 1996.
\bibitem{Da Prato2} G. Da Prato, \emph{Introduction to stochastic analysis and Malliavin calculus}, Edizioni Della Normale, 2014.
%\bibitem{Luca_Li} L. Dieci, W. Li and H. Zhou, \emph{ A new model for realistic perturbation of stochastic oscillators}, Journal of Differential Equations, 261 (2016) 2502 - 2527.
\bibitem{Duan} J. Duan, K. Lu and B. Schmalfuss, \emph{Invariant Manifolds for Stochastic Partial Differential Equations, } The Analysis of Probability, Vol. 31 (2003), 2109 - 2135.
%\bibitem{Dons76} M. D. Donsker and S. R. S Varadhan, \emph{Asymptotic Evaluation of Certain Markov Process Expectations for Large Time -- III}, Communications on Pure and Applied Mathematics, Vol. XXIX (1976), 389 -- 461.
%\bibitem{Dupuis07} P. Dupuis and R. S. Ellis, \emph{A weak convergence approach to the theory of large deviations,} Wiley-Interscience, 1997.
%\bibitem{Dupuis16} P. Dupuis, M. A. Katsoulakis, Y. Pantazis and P. Plech\'ac, \emph{Path-space information bounds for uncertainity quantification and sensitivity analysis of stochastic dynamics,} SIAM Journal of Uncertainty Quantification, Vol. 4 (2016), 80 -- 111.
\bibitem{Elw1} K. D. Elworthy, \emph{Stochastic dynamical systems and their flows.} In: Stochastic Analysis, ed. A. Friedman, M. Pinsky, London-New York press, (1978) 79 - 95.
\bibitem{ElwF} K. D. Elworthy and X.-M. Li, \emph{Formulae for the derivative of heat semigroups}, Journal of functional analysis, Vol. 125 (1994), 252 - 286.
\bibitem{Chung} C. Feng, H. Zhao and B. Zhou, \emph{Pathwise random periodic solutions of stochastic differential equations,} Journal of Differential Equations, Vol. 251 (2011), 119 - 149.
\bibitem{Feng} C. Feng and H. Zhao, \emph{Random Periodic Solutions of SPDEs via Integral Equations and Wiener-Sobolev Compact Embedding,} Journal of Functional Analysis, Vol. 262 (2012)  4377 - 4422.
%\bibitem{Zhao3} C. Feng, Y. Liu and H. Zhao, \emph{Numerical approximation of random periodic solutions of stochastic differential equations}, Z. Angew. Math. Phys. (2017) 68: 119.
\bibitem{Zhao2} C. Feng and H. Zhao, \emph{Random periodic processes, periodic measures and ergodicity,} available at arXiv:1408.1897v4 [math.PR] 
%\bibitem{Figali} A. Figali, \emph{Existence and uniqueness of martingale solutions for SDEs with rough or degenerate coefficients,} Journal of Functional Analysis, 254 (2008) 109 -- 153.
%\bibitem{Flandol3} F. Flandoli, B. Gess and Michael Scheutzow, \emph{Synchronisation by noise,} Probability Theory and Related Fields, 168 (2017), 511 -- 556.
%\bibitem{Garrido03} M. I. Garrido and J. A. Jaramillo, \emph{Homomorphisms on functin lattices}, Monatsh. Math. 141 (2004), 127 -- 146.
%\bibitem{Garrido08} M. I. Garrido and J. A. Jaramillo, \emph{Lipschitz-type functions on metric spaces}, Journal of Mathematical Analysis and Applications, 340 (2008), 282 -- 290.
%\bibitem{Froyland17} G. Froyland and P. Koltai, \emph{Estmating long-term behaviour of periodically driven flows without trajectory integration}, Nonlinearity 30 (2017) 1948 -- 1986.
%\bibitem{Giesl} P. Giesl and S. Hafstein, \emph{Local Lyapunov functions for Periodic and Finite-Time ODEs,} In: Recent Trends in Dynamical Systems,  Proceedings of a Conference in Honor of J\"urgen Scheurle, Springer 2013.
\bibitem{Har1} T. E. Haris, \emph{Brownian motions on the homeomorphisms of the plane,} Annals of Probability, Vol. 9 (1981), 232--254.
\bibitem{Hasm} R. Z. Has\'minskii, \emph{Stochastic stability of differential equations,} Sijthoff and Noordhoff, 1980 (Translated from Russia).
\bibitem{Ikeda} N. Ikeda and S. Watanabe, \emph{Stochastic Differential Equations and Diffusion Processes,} North Holland - Kodansha, Tokyo, 1981. 
\bibitem{Kak} S. Kakutani, \emph{Random Ergodic Theorems and Markov Processes with a Stable Distributions,} Proceedings of second Berkeley Symposium on Mathematical Statistics and Probability (Univ. of Calif. press, 1951),  247 - 261.
\bibitem{Hkunita} H. Kunita, \emph{Stochastic differential equations and stochastic flow of diffeomorphisms,} \'Ecole d'\'et\'e de Probabilit\'es de Saint-Flour 12, 1982 Lecture Notes in Mathematics, Vol. 1097 (1984), 143 - 303. 
\bibitem{Kunita} H. Kunita, \emph{Stochastic flows and stochastic differential equations,} Cambridge University press, Cambridge, 1990.
\bibitem{Kurr} C. Kurrer and K. Schulten, \emph{Effect of noise and perturbations of limit cycle systems,} Physica D, Vol. 50 (1991), 311 - 320.
%\bibitem{Kushn} H. Kushner, \emph{Stochastic stability and control,} Academic press, New York, 1967.
\bibitem{Lien}  A. Liénard, \emph{Etude des oscillations entretenues,} Rev. Gén. d’Elect. 23 (1928), 901 - 912 and 946 - 956.
%\bibitem{Igal} F. Liese and I. Vadja, \emph{On Divergences and Informations in Statistics and Information Theory,} Transcations on Information Theory, Vol. 52 (2006), 4394 -- 4412.
%\bibitem{Luo} D. Luo, \emph{Quasi-invariance of Lebesgue measure under the homeomorphic flow generated by SDE with non-Lipschitz coefficient,} Bull. Sc. Math. 133 (2009) 205 -- 228.
\bibitem{Mao} X. Mao, \emph{Exponential stability of stochastic differential equations}, Marcel Dekker, New York, 1994.
\bibitem{Mattingly} J. C. Mattingly, A. M. Stuart and D. J. Higham, \emph{Ergodicity for SDEs and approximations: locally Lipschitz vector fields  and degenerate noise,} Stochastic Processes and their Applications, Vol. 101 (2002), 185 - 232.
\bibitem{Meyn93} S. P. Meyn and R. L. Tweedie, \emph{Stability of Markovian processes III:Foster-Lyapunov criteria for continuous-time processes}, Advances in Applied Probability, Vol. 25 (1993), 518 - 548.
\bibitem{Meyn931} S. P. Meyn and R. L. Tweedie, \emph{Stability of Markovian processes II: Continuous-time processes and sampled chains,} Advances in Applied Probability, Vol. 25 (1993), 487 -517.
%\bibitem{Meyn92} S. P. Meyn and R. L. Tweedie, \emph{Stability of Markovian processe I: Criteria for discrete-time chains}, Advances in Applied Probability, Vol. 24 (1992), 542 -- 574.
\bibitem{Ponc} H. Poincar\'e, \emph{Memoire sur les courbes definier par une equation differentiate}, Jour. Math. Pures et appli. 3, tome 7 (1881), 375-422; 3, tome (1882), 251-296; 4, tome 1 (1885), 167-244; 4, tome 2 (1886), 151-217.
\bibitem{Schmalfuss} B. Schmalfuss, \emph{Lyapunov functions and non-trivial stationary solutions of stochastic differential equations,} Dynamical Systems, Vol. (2001), 303 - 317.
%\bibitem{Schurz} H. Schurz, \emph{Verification of Lyapunov functions for the analysis of stochastic Li\'enard equations,} Journal of Sound and Vibration, Vol. 325 (2009) 938-949.
%\bibitem{Rsmith} R. A. Smith, \emph{Orbital Stability for Ordinary Differential Equations,} Journal of Differential Equations 69, 265 - 287 (1987).
%\bibitem{Stroock} D. W. Stroock and S. R. S Varahdan, \emph{Multidimensional Diffusion Processes,} Springer-Verlag, 1979.
\bibitem{Uda} K. Uda, \emph{Existence of random invariant periodic curves via random semiuniform ergodic theorem,} Stochastics and Dynamics, Vol. 17 (2017), 175007.
\bibitem{Ulam} S. M. Ulam and J. Von Neumann, \emph{Random Ergodic Theorems,} Bulletin of American Mathematical Society, Vol. 51 (1945), 660.
\bibitem{Van} B. van der Pol, \emph{On Relaxation Oscillations,} Philosophical Magazine Series 7, Vol. 2 (1926), 978 - 992.
%\bibitem{Kampen} N. G. Van Kampen, \emph{Stochastic Processes in Physics and Chemistry,} Elsevier, 2007.
%\bibitem{Wan1} B. Wang, \emph{Periodic and almost periodic random inertial manifolds for non-autonomous stochastic equations,}  page 189 - 208, Continuous and Distributed Systems II: Theory and Applications, Edited by V. A. Sadovnichiy and M. Z. Zgurovsky, Springer, 2015.
 \bibitem{Wan2} B. Wang, \emph{Stochastic bifurcation of pathwise almost random periodic and almost automorphic solutions for random dynamical systems,} Discrete and Continuous Dynamical Systems, Vol. 36 (2015), 3745 - 3769.
 \bibitem{Wan3} B. Wang, \emph{ Existence, stability and bifurcation of random complete and periodic solutions of stochastic parabolic equations,} Nonlinear Analysis: Theory, Methods and Applications, Vol. 103 (2014), 9 - 25.
\bibitem{Jeffrey} J. B. Weiss and E. Knobloch, \emph{A Stochastic Return Map for Stochastic Differential Equations,} Journal of Statistical Physics, Vol. 58 (1990), 863 - 883.
\bibitem{Huaiz} Q. Zhang and H. Zhao, \emph{ Stationary Solutions of SPDEs and infinite horizons BDSDEs,} Journal of Functional Analysis, Vol. 252 (2007), 171 - 219. 
\bibitem{Zhao} H. Zhao and Z. Zheng, \emph{Random Periodic Solutions of Random Dynamical systems,} Journal of differential equations, Vol. 246 (2009), 2020 - 2038. 
 
\end{thebibliography}
\end{document}